\documentclass[a4, 11pt]{article}

\usepackage[english]{babel}
\usepackage{graphicx}
\usepackage{xcolor}
\usepackage{amsmath,amsfonts,amssymb,amsthm}
\usepackage[linesnumbered,ruled,vlined]{algorithm2e}
\SetKwInOut{Repeat}{Repeat}

\usepackage{hyperref}

\usepackage{url}
\usepackage{cleveref}
\usepackage{caption}
\usepackage{subcaption}
\usepackage{graphicx}

\newcommand{\C}{\mathbb{C}}
\newcommand{\R}{\mathbb{R}}
\newcommand{\Rp}{\mathbb{R}_{\ge 0}}
\newcommand{\N}{\mathbb{N}}

\newcommand{\cA}{\mathcal{A}}
\newcommand{\cC}{\mathcal{C}}
\newcommand{\cD}{\mathcal{D}}
\newcommand{\cF}{\mathcal{F}}
\newcommand{\cG}{\mathcal{G}}
\newcommand{\cL}{\mathcal{L}}
\newcommand{\cN}{\mathcal{N}}
\newcommand{\cS}{\mathcal{S}}
\newcommand{\cU}{\mathcal{U}}
\newcommand{\cX}{\mathcal{X}}

\newcommand{\ve}{\varepsilon}

\def\MIN#1{\min_{#1} \ }
\def\INF#1{\inf_{#1} \ }
\def\MAX#1{\max_{#1} \ }
\def\argmin#1{\underset{#1}{\arg \min} \ }

\DeclareMathOperator{\diag}{\rm diag}
\newcommand{\setdef}[2]{\left\lbrace\ #1\ \left|\ \vphantom{#1} #2\ \right.\right\rbrace}

\begin{document}

\title{Distributed Optimization for Energy Grids: \\ A Tutorial on ADMM and ALADIN}

\author{Lukas Lanza\thanks{Optimization-based Control Group, Institute of Mathematics, Technische Universität Ilmenau, Ilmenau, Germany. \textsc{E-Mail}: \url{lukas.lanza@tu-ilmenau.de}} 
\and 
Timm Faulwasser\thanks{Present address: Institute of Control Systems, TU Hamburg, Hamburg, Germany. \textsc{E-Mail}: \url{timm.faulwasser@ieee.org}. This work has been conducted while TF was with the Institute of Energy Systems, Energy Efficiency and Energy Economics, TU Dortmund University, Dortmund, Germany.} 
\and 
Karl Worthmann\thanks{Optimization-based Control Group, Institute of Mathematics, Technische Universität Ilmenau, Ilmenau, Germany. \textsc{E-Mail}: \url{karl.worthmann@tu-ilmenau.de}}
}
\maketitle

\textbf{Abstract}: The ongoing transition towards energy and power systems dominated by a large number of renewable power injections to the distribution grid poses substantial challenges for system operation, coordination, and control. Optimization-based methods for coordination and control are of substantial research interest in this context. Hence, this chapter provides a tutorial introduction of distributed optimization algorithms for energy systems with a large share of renewables. Specifically, we focus on the Alternating Direction Method of Multipliers (ADMM) and on the Augmented Lagrangian Alternating Direction Inexact Newton (ALADIN) method as both algorithms are frequently considered for coordination and control of power and energy systems. We discuss the application of ALADIN and ADMM to AC optimal power flow problems and to energy management problems. Moreover, we give an outlook on open problems.

\section{Introduction}
The ongoing energy transition towards power and energy systems involving a large number of distributed resources and a rapidly growing share of volatile renewable energy production induces multiple challenges for operation, coordination, and control of such systems. 
Traditionally the flow of electrical power  is mainly unidirectional, i.e., from the large fossil power plants feeding into transmission systems to the customers supplied by the distribution grids.  Networks with increasing numbers of distributed resources (prosumers, battery storage and controllable renewable generation) exhibit bidirectional flow of power flow and the established coordination mechanisms between transmission and distribution systems have to be rethought, see, e.g., \cite{westermann2019curative,boie2016opportunities} and  references therein.

This chapter discusses solution concepts for the situation described above from the mathematical point of view.
In particular, we focus on coordination and control via numerical optimization. We consider \emph{distributed optimization algorithms} as these methods allow to tackle the problem of coordinating a large number of individual systems while avoiding  the clustering of all information in a single centralized entity. Specifically, we discuss the Alternating Direction Method of Multipliers (ADMM) and  the Augmented Lagrangian Alternating Direction Inexact Newton (ALADIN) method as both are frequently considered for power and energy systems.  

An alternative route would be complexity reduction via model reduction, see, e.g., \cite{rasheduzzaman2015reduced,gnarig2022model}. However, this typically results in sub-optimal performance. Moreover, in the context of frequency and voltage stabilization also the design of distributed sub-optimal feedback laws is considered~\cite{yazdanian2014distributed,zhao2015distributed,ortega2018frequency}.

The remainder of this chapter is structured as follows:
In Section \ref{Sec:dOPT} we explain the conceptual idea of distributed optimization including a brief exposition of ADMM and ALADIN. Section \ref{Sec:Applications} focuses on two illustrative use cases for distributed optimization in energy systems: energy management and AC Optimal Power Flow (AC-OPF). This chapter ends with an outlook on open problems in Section \ref{Sec:Outlook}.

\subsubsection*{Notation}
Within this chapter we use
$\langle \cdot, \cdot \rangle$ as the scalar product on~$\R^n$; 
$\| \cdot \|_p$ denotes the p-norm in~$\R^n$, and for brevity, we use $\| \cdot\| := \|\cdot\|_2$.
Moreover, $\|\cdot\|^2_Q := \langle \cdot , Q \cdot \rangle$ for a symmetric positive-definite matrix~$Q$.
Further, $I_n := \diag(1,\ldots,1) \in \R^{n \times n}$ denotes the identity matrix of dimension $n \in \N$.
In the optimization problems, we write the Lagrangian multiplier after the constraints, i.e., ``S~$| \lambda$'' means that~$\lambda$ is the Lagrange multiplier for the constraint~S.
For $p\in\N $ the space of $p,m,n$-times continuously differentiable
functions on $V \subset \R^m$ with image in $\R^n$ is denoted by $C^p(V,\R^n)$.

\section{Distributed optimization} \label{Sec:dOPT}
Optimization in general aims at finding solutions which are superior to other solutions in a particular sense.
The notion of optimality is defined by asking a preassigned performance criterion to be minimal (or maximal).
In the context of energy systems and power supply, typical optimality criteria  
are related to 
transmission losses~\cite{murray2018hierarchical,Sauerteig2021},
generation costs~\cite{faulwasser2018optimal}, 
efficiency~\cite{parisio2014model}, incorporation of small-scale renewables and storage devices~\cite{muhlpfordt2019chance,braun2020towards,hu2023economic}, building control~\cite{oldewurtel2012use}, or voltage stability~\cite{schiffer2015voltage,bai2022voltage}, see also~\cite{dorfler2015breaking} as well as the surveys~\cite{schiffer2016survey,molzahn2017survey} and references therein. 
For example, the power grid should be operated such that  
line losses are minimized, at the same time as much  
energy from residential photovoltaic generation is used while local storage devices are limited by capacity and charging constraints.
Therefore, the resulting optimization problems 
are of large size and of high complexity with many individual cost functions and constraints.
To solve these problems efficiently, it is necessary to investigate the structure and to consider, e.g., whether the problems can be decomposed into subproblems to be solved in parallel.
For the sake of clarity, we briefly state our notion of decentralized and distributed optimization.
Following~\cite{stomberg2023cooperative}, in decentralized optimization, the overall problem can be split into subproblems, which are coupled such that each subproblem can be solved separately
while only communicating its solution to the \textit{coupled} subproblems, see, e.g., Problem~\eqref{eq:ToyProblem1} below.
The subproblems communicate with each other directly, in particular, no extra communication entity is involved.
In distributed optimization, the overall problem can be split
into subproblems, which are coupled with each other via constraints or via the cost functions, see e.g., Problem~\eqref{eq:ToyProblem2} below. 
Each subproblem subject to the connecting constraints is solved individually. A central entity, e.g., a market maker~\cite{worthmann2014distributed,faulwasser2018optimal}, communicates the local solutions and solves, e.g., a consensus optimization problem based on the local solutions.
With reference to~\cite[Sec.~I]{molzahn2017survey}, we emphasis that the above terms are used ambiguously.

We highlight two main aspects leading to distribution of optimization problems.
First, the computational tasks to be done can be distributed to several units.
This means that a large-scale optimization problem is split into several small-scale subproblems, where
each subproblem is amenable to optimization algorithms individually.
Therefore, the computation task can be done in parallel, which results in less computation time.
The parallelization may even be necessary to make the problem solvable at all, as the overall problem would have required too much memory, for example.
A concise overview of distributed optimization, with regard to energy systems, can be found in, e.g., \cite{molzahn2017survey,yang2019survey}.
Since the overall optimization problem typically does not consist of completely decoupled subproblems,
the artificially split subproblems have to be equipped with coupling constraints to obtain an equivalent problem, cf. the ideas presented in, e.g.,~\cite{rantzer2009dynamic,chatzipanagiotis2015augmented,houska2016augmented,falsone2017dual}.
This idea is discussed in \Cref{Sec:Decomposition}. \
The second main aspect of distributing the optimization problem is linked to data security and sovereignty.
Solving the overall problem as a whole, the computation entity must have access to all constraints, all objectives, and all states; in particular, all local information has to be shared.
If, however, the overall problem is split into several subproblems, which can be solved individually in parallel, only little information, e.g., local optimal values, has to be shared with a central entity since the latter only has to solve a consensus problem, see, for example, the works~\cite{molzahn2017survey,varagnolo2015newton,wang2023differentially}. 
This means that privacy is maintained to a much greater extent. 

Another possible option to take care of data security is to encrypt local information, share it with the computation entity which operates on the encrypted data, receive a still encrypted solution, locally decrypt the data and apply the solution, cf.~\cite{darup2021encrypted,schluter2023brief}. Although these techniques are subject of ongoing research, they are still quite limited in terms of the complexity of the computational tasks.

Besides parallel computation and higher level of data privacy, a distributed optimization problem should
have the following properties.
The problem allows plug-and-play, i.e., including new components (e.g., PV or storage devices) does not change the
overall problem, but only few modifications are required, cf.~\cite{mbuwir2020distributed,zeilinger2013plug,dorfler2014plug,crisostomi2014plug}. 
Moreover, the solution techniques are scalable in the sense that the same algorithms are applicable if the problem is extended, e.g. by medium-scale energy generation.
Furthermore, simple manipulations such as utilizing averaged quantities allow to keep the number of unknowns~\cite{braun2018hierarchical} and the communication effort~\cite{jiang2020distributed} independent of the number of subsystems.

\ \\
Distributed optimization aims at separating an large-scale  and/or structured problem into
several small-scale subproblems connected via coupling constraints.
This  typically leads to partially separable problems of the form
\begin{equation} \label{eq:SeparableCoupledProblem}
    \begin{aligned}
        \MIN{x} \sum_{i=1}^N k_i(x_i) \ \ \text{subject to} \ \ 
        \begin{cases}
            \sum_{i=1}^N A_i x_i = b, \\
             h_i(x_i) \le 0, & i \in \{1,\ldots,N \},
        \end{cases}
    \end{aligned}
\end{equation}
where~$k_i: \R^n \to \R$, $h_i: \R^n \to \R^{n_h}$, $i \in \{1,\ldots,N\}$ are to be specified in the respective context, e.g., their regularity and convexity properties.
The matrices~$A_i \in \R^{m \times n}$ and the vector~$b \in \R^m$ couple the subproblems.

\ \\
Before we present two approaches to solve optimization problems~\eqref{eq:SeparableCoupledProblem}
in \Cref{Sec:ADMMandALADIN}, namely the algorithms ADMM and ALADIN, we first illustrate simple and briefly the idea of
decomposing optimization problems.

\subsection{Decomposition techniques} \label{Sec:Decomposition}
In this section, we briefly discuss two decomposition strategies, namely decomposition in primal variables, and dual decomposition.
The underlying idea is to split a complex optimization problem into more simple subproblems and a connecting master problem.
Roughly speaking the main difference is that decomposition in primal variables leads to subproblems where the connecting variables are primal variables, while dual decomposition yields subproblems with duals as connecting variables.
The difference can be illustrated as follows:
Consider the case where usage of several resources is coupled in a complex way.
Decomposition in primal variables yields a set of subproblems where the connecting master problem directly controls the resources.
In contrast, via dual decomposition the master problem controls the resources' prices, which correspond to the Lagrange multipliers.

\paragraph{Decomposition in primal variables}
Consider a structured optimization problem
\begin{equation} \label{eq:ToyProblem1}
    \MIN{x,y,z} K(x,y,z) = k_1(x,z) + k_2(y,z),
\end{equation}
where~$x,y$ are local variables, and~$z$ are coupling variables.
To solve~\eqref{eq:ToyProblem1}, fix~$z$ and simultaneously solve the subproblems 
\begin{align*}
    \MIN{x} k_1(x,z) \quad \text{and} \quad \MIN{y} k_2(y,z),
\end{align*}
and obtain optimal $k_1^\star(z)$, $k_2^\star(z)$.
Problem~\eqref{eq:ToyProblem1} is then equivalent to the \textit{primal master problem}
\begin{equation*}
    \MIN{z} k_1^\star(z) + k_2^\star(z),
\end{equation*}
with coupling variable~$z$.
Each iteration of solving the primal master problem requires solving the subproblems.
Note that, if the original problem~\eqref{eq:ToyProblem1} is convex, then the primal master problem is convex.

The decomposition in primal variables is reminiscent of Bellman's Dynamic Programming. Indeed, structured decomposition of path and tree-structured problems can be achieved via decomposition in primal variables and recursion (which alleviates the need for a master problem), see, e.g., \cite{jiang2021decentralized} and \cite{engelmann2022scalable,engelmann2023approximate} and references therein.

\paragraph{Dual decomposition}
If the optimization problem is subject to further constraints, the decomposition can be applied to the dual problem.
Originally proposed in~\cite{everett1963generalized}, the idea is to include the constraints via Lagrange multipliers and to solve a dual ascent problem.

Consider
\begin{equation} \label{eq:ToyProblem2}
    \begin{aligned}
        \min_{x,y,z_1,z_2} K(x,y,z_1,z_2) = k_1(x,z_1) + k_2(y,z_2), \ \ \text{subject to} \ \ z_1 = z_2,
    \end{aligned}
\end{equation}
with local variables~$x,y$, and coupling variables~$z_1, z_2$ to satisfy so-called consensus.
We define the corresponding Lagrangian function
\begin{equation*} 
    L(x,y,z_1,z_2,\lambda) := K(x,y,z_1,z_2) + \lambda(z_1-z_2),
\end{equation*}
where~$\lambda$ is the so-called Lagrange multiplier.
Denoting $p^\star := \INF{x,y,z} k_1(x,z) + k_2(y,z)$ we may estimate with the dual function
\begin{equation*}
    D(\lambda) := \MIN{x,z} (k_1(x,z) - \lambda z) + \MIN{y,z} ( k_2(y,z) + \lambda z) \le p^\star 
\end{equation*}
a lower bound for~$p^\star$ for all~$\lambda$.
Then, the dual problem
\begin{equation*}
    \MAX{\lambda} D(\lambda) = \MAX{\lambda} \left\{ \MIN{x,z}(k_1(x,z) + \lambda z) + \MIN{y,z}( k_2(y,z) -\lambda z)  \right\}
\end{equation*}
provides the tightest lower bound for~$p^\star$.
Note that if the primal problem~\eqref{eq:ToyProblem2} is convex and feasible, then $D(\lambda^\star) = p^\star$.
However, usually in practice the so-called duality gap $D(\lambda^\star) < p^\star$ exists, even for optimal dual~$\lambda^\star$, that is, sub-optimality of solutions.
Moreover, if $z_1^\star \neq z_2^\star$, then it is not straightforward to recover the primal optimal $(x^\star , y^\star, z_1^\star, z_2^\star)$.

\ \\
The method of dual decomposition has been applied to various optimization problems in systems and control, see, e.g.,~\cite{rantzer2009dynamic,giselsson2013accelerated,falsone2017dual,botkin2021distributed}.

\subsection{ADMM and ALADIN} \label{Sec:ADMMandALADIN}
In this section we present and discuss two algorithms, namely ADMM and ALADIN, to solve distributed optimization problems~\eqref{eq:SeparableCoupledProblem}, 
where~$k_i, h_i$ are local costs and constraints, i.e., dependence is only on the local variables~$x_i$.
The matrices~$A_i \in \R^{m \times n}$ and the vector~$b \in \R^m$ are used to model the coupling between subsystems.
Regularity and convexity requirements of~$k_i, h_i$ for each algorithm will be discussed in the next sections. 

\paragraph{ADMM}
The acronym ADMM refers to the Alternating Direction Method of Multipliers, 
originally proposed in~\cite{glowinski1975approximation,gabay1976dual}.
It is a further development of dual decomposition.
The idea can be illustrated as follows.
For $k_i : \cD \to \R\cup\{+\infty\}$, $i=1,2$ consider the simplified problem 
with consensus constraints
\begin{equation} \label{eq:ToyProblem3}
    \begin{aligned}
        \min_{z_1,z_2 \in \cD} K(z_1,z_2) = k_1(z_1) + k_2(z_2), \ \ \text{subject to} 
       \ \ z_1 = z_2.
    \end{aligned}
\end{equation}
This problem is solved approximately by alternating solve for~$z_1$ with~$z_2$ fixed, and then solve for~$z_2$ with~$z_1$ fixed.
For a regularization (augmentation) parameter~$\mu > 0$ we define the augmented Lagrangian
\begin{equation} \label{eq:ToyLagrangian}
    L_\mu(z_1,z_2,\lambda) := k_1(z_1) + k_2(z_2) + 
    \langle z_1 - z_2, \lambda \rangle + \tfrac{\mu}{2} \| z_1 - z_2  \|^2 .
\end{equation}
Then the ADMM algorithm essentially repeats the following steps
\begin{equation*} 
\begin{aligned}
    z_1^k & = \argmin{z_1} L(z_1,z_2^{k-1},\lambda^{k-1}) , \\
    z_2^k & = \argmin{z_2} L(z_1^k,z_2,\lambda^{k-1}) , \\
    \lambda^k & = \lambda^{k-1} + \mu(z_1^k - z_2^k),
\end{aligned}
\end{equation*}
and so approaches the optimal solution, see also~\cite[Alg.~2]{stomberg2022decentralized} for a concise presentation of the algorithm.
To solve the more general problem~\eqref{eq:SeparableCoupledProblem}, we recall in \Cref{algo:ADMM} below the version of ADMM from~\cite{houska2016augmented}.
\begin{algorithm} [ht!]
\KwIn{Initial guess~$x_i \in \R^n$ and~$\lambda_i \in \R^m$, penalty parameter~$\rho > 0$, and numerical tolerance~$\ve > 0$}
\Repeat{}
 \label{ADMM:Step1}
        Solve for $i \in \{1,\ldots,N\}$ the decoupled problems
        \begin{equation*} 
        \MIN{y_i} k_i(y_i) + \langle A_i y_i, \lambda \rangle + \tfrac{\rho}{2} \| A( y_i - x_i) \|^2 \ \ \text{subject to} \ \ h_i(y_i) \le 0 . \tag{NLP}
        \end{equation*} 
    \\ 
    If $\| \sum_{i=1}^N A_i y_i - b \|_1 \le \ve$, then terminate with $x^\star = y$.
    \\
    Implement the dual gradient steps $\lambda^+_i = \lambda_i + \rho A_i(y_i - x_i)$.
    \\
    Solve the coupled QP \label{ADMM:Step4}
    \begin{equation*}  
        \MIN{x^+} \sum_{i=1}^N \frac{\rho}{2} \| A_i (y_i -x_i^+) \|^2 + \langle A_i x_i^+, \lambda_i^+ \rangle 
        \quad \text{subject to} \quad \sum_{i=1}^N A_i x_i^+ = b. \tag{QP}
    \end{equation*}
    \\
    Update the iterates $x \leftarrow x^+$ and $\lambda \leftarrow \lambda^+$, and go to Step~1. \\
\KwOut{Numerical optimal solution~$x^\star$}
\caption{ADMM in virtue of~\cite{houska2016augmented}}
\label{algo:ADMM}
\end{algorithm}
In~\cite[Sec.~3.2, App.~A]{boyd2011distributed} it has been proven for separable problems~\eqref{eq:ToyProblem3} with  closed, proper, and convex objective functions~$k_1,k_2$, under a saddle-point condition on the un-augmented Lagrangian ($\mu = 0$), that
the iterates are feasible ($(z_1^n-z_2^n)_{n \in \N} \to 0$ as~$n \to \infty$), 
the costs converge to an optimal value ($(K(z_1^n,z_2^n,\lambda^n))_{n \in \N} \to K^\star$ as~$n \to \infty$),
and that the duals are optimal too ($\lambda^n \to \lambda^\star$ as~$n \to \infty$).
ADMM (\Cref{algo:ADMM}) has been successfully applied to optimization problems arising in energy systems, see,  e.g.,~\cite{erseghe2014distributed,magnusson2015distributed,braun2018hierarchical}. Moreover, one should note that ADMM as such is rather a family of algorithms than a single rule, cf. the different expositions and variants in \cite{boyd2011distributed,stomberg2022compendium, braun2016distributed}.

\paragraph{ALADIN}
The acronym ALADIN stands for Augmented Lagrangian Alternating Direction Inexact Newton.
First proposed in~\cite{houska2016augmented},
this algorithm is designed to solve non-convex problems~\eqref{eq:SeparableCoupledProblem},  
where we assume regularity of the costs $k_i \in C^2(\R^n; \R)$ as well as for the constraints $h_i \in C^2(\R^n; \R^{n_h})$ for all $i \in \{1,\ldots,N\}$.
Furthermore, we assume that~\eqref{eq:SeparableCoupledProblem} is feasible and all local minimizers are regular KKT points.
Note that, in contrast to ADMM, convexity of the functions~$k_i,h_i$ is not required.
Importantly, the functions~$k_i, h_i$ take only~$x_i$, this means, the problem is separable into~$N$ subproblems, which are coupled via the consensus constraint~$[A_1, \ldots, A_N] x = b$.
Similar to ADMM, the idea is to decompose~\eqref{eq:SeparableCoupledProblem} into subproblems, which can be solved in parallel, where the coupling constraints are formulated as costs.
The first two steps of \Cref{algo:ALADIN} (ALADIN) solve the decoupled subproblems in primal variables analogous to \Cref{algo:ADMM} (ADMM).
In Step~3 an approximation of the constraint Jacobian is obtained. 
The latter is similar to inexact SQP methods and it allows to cope with inexact constraint Jacobian and resulting errors in the step direction.
Moreover, in contrast to ADMM, an approximation of the Hessian of $k_i(x_i) + \langle h_i(x_i), \kappa_i \rangle$ is required to solve the coupling problem~\eqref{eq:dQP} in Step~4.
For the update of~$x$ and~$\lambda$ in Step~5 the line search parameters can be set~$\alpha_1=\alpha_2=\alpha_3=1$, or can be computed according to~\cite[Alg.~3]{houska2016augmented} to obtain global convergence guarantees.
\Cref{algo:ALADIN} summarizes the ALADIN as proposed by~\cite{houska2016augmented}.

\begin{algorithm}[ht!] 
\KwIn{Initial guess for~$x_i \in \R^n$ and~$\lambda \in \R^m$, and numerical tolerance~$\ve > 0$}
\Repeat{}
 \label{ALA:Step1}
        Choose penalty parameter~$\rho \ge 0$ and $0 \le \Sigma_i \in \R^{n \times n}$, 
        and solve for $i \in \{1,\ldots,N\}$ the decoupled problems
        \begin{equation*} \label{eq:dOP}
        \MIN{y_i} k_i(y_i) + \langle A_i y_i, \lambda \rangle + \tfrac{\rho}{2} \| y_i - x_i \|^2_{\Sigma_i}, \ \ \text{subject to} \ \ h_i(y_i) \le 0 \ \ \ | \, \kappa_i  \tag{NLP}
        \end{equation*} 
        to either local or global optimality.
    \\ 
    If $\| \sum_{i=1}^N A_i y_i - b \|_1 \le \ve$ and $\| \Sigma_i (y_i -x_i)\|_1 \le \ve $, then terminate with $x^\star = y$.
    \\
    Choose for $i \in \{ 1,\ldots,N \}$ an approximation of the Jacobian $C_i \approx C_i^\star $ 
    \begin{equation*}
    \forall \, j = 1, \ldots, n_h  \, : \ \ 
    C_{ij}^\star = 
    \begin{cases}
    \nabla_x (h_i(x))_j |_{x = y_i}, & \text{if} \ (h_i(y_i))_j = 0, \\
    0, & \text{otherwise},
    \end{cases}
    \end{equation*}
    and compute $g_i = \nabla k_i(y_i) + \langle C_i^\star - C_i),  \kappa_i \rangle$, and choose a symmetric positive-definite approximation of the Hessian \label{ALA:Step3}
    \begin{equation*}
    H_i \approx \nabla^2 \big( k_i(x_i) + \langle h_i(y_i), \kappa_i \rangle \big) \in \R^{n \times n}.
    \end{equation*}
    \\
    For $\mu > 0$ solve the coupled QP \label{ALA:Step4}
    \begin{equation*} \label{eq:dQP}
    \begin{aligned}
        & \MIN{\Delta y, s} \sum_{i=1}^N  \tfrac{1}{2} 
        \| \Delta y_i\|^2_{H_i}
        + \langle \Delta y_i, g_i \rangle   + \langle s, \lambda \rangle + \tfrac{\mu}{2} \| s\|^2_2  \\
        & \text{subject to} \ \
        \begin{cases}
            \sum_{i=1}^N A_i (y_i + \Delta y_i) = b +s \ \ \ | \,  \lambda_{\rm QP}, \\
            C_i \Delta y_i = 0, \quad i \in \{1,\ldots,N\} . 
        \end{cases}
    \end{aligned}
    \tag{QP}
    \end{equation*}
    \\
    For step-sizes $\alpha_1, \alpha_2, \alpha_3 \ge 0$ set
    \begin{equation*}
    x^+ = x + \alpha_1 (y-x) + \alpha_2 \Delta y, \quad \lambda^+ = \lambda + \alpha_3 (\lambda_{\rm QP} - \lambda) .
    \end{equation*}
    \\
    Update the iterates $x \leftarrow x^+$ and $\lambda \leftarrow \lambda^+$, and go to Step~1. \\
\KwOut{Numerical optimal solution~$x^\star$}
\caption{ALADIN~\cite{houska2016augmented}}
\label{algo:ALADIN}
\end{algorithm}

\paragraph{Comparing ADMM and ALADIN}
In both algorithms, \Cref{algo:ADMM} (ADMM) and \Cref{algo:ALADIN} (ALADIN), the~$N$ subproblems~\eqref{eq:dOP} are nonlinear, but compared to the overall problem of small scale. The coupling problems are of large scale, but only linear equality constraints are involved.
Moreover, due to the separated structure of~\eqref{eq:dOP}, no detailed information about the costs~$k_i$, the constraints~$h_i$ and the coupling~$A_i$ have to be shared to the central entity.
For~$\rho = 0$ the method in \Cref{algo:ALADIN} is equivalent to dual decomposition, see \cite[App.~A]{houska2016augmented} for detailed analysis.
For $H_i = \rho A_i^\top A_i$, $\Sigma_i = A_i^\top A_i$, $C_i = 0$, and large values of~$\mu$ ($\mu \to \infty$) similarity between \Cref{algo:ADMM} and \Cref{algo:ALADIN} is in detail discussed in~\cite[Sec.~5]{houska2016augmented}.

Unlike ADMM, ALADIN can cope with non-convex costs~$k_i$ as well as with non-convex constraints~$h_i$, cf.~\cite{houska2016augmented}.
Invoking an $L_1$-penalty function, it is possible to choose~$\alpha_1,\alpha_2,\alpha_3$, cf.~\cite[Algo.~3]{houska2016augmented}, such that \Cref{algo:ALADIN} converges globally to local minimizers of~\eqref{eq:SeparableCoupledProblem}, cf.~\cite[Lem.~3]{houska2016augmented}.
It turns out that in each iteration either $\alpha_1=\alpha_2 = \alpha_3 = 1$, or $\alpha_1 = \alpha_2 = 0$ and~$\alpha_3 = 1$, cf.~\cite[Thm.~4]{houska2016augmented}.
If~\eqref{eq:SeparableCoupledProblem} is feasible and bounded from below, $\alpha_1,\alpha_2 , \alpha_3$ computed by~\cite[Algo.~3]{houska2016augmented}, and~$\rho$ large enough, then \Cref{algo:ALADIN} terminates after finite steps, cf.~\cite[Thm.~2]{houska2016augmented}.
As elaborated in~\cite[Sec.~III]{engelmann2018toward} the approximation of the Hessian in Step~\ref{ALA:Step3} can be replaced by the damped Boyden-Fletcher-Goldfarb-Shanno (BFGS) update, which only involves previous gradient updates of the Lagrangian.
Using this update also reduces the overall communication effort.

In \Cref{Tab:CompareADMMandALADIN} we summarize a qualitative comparison of ADMM and ALADIN algorithms.
\begin{table}[ht!]
    \caption{
    Comparison of \Cref{algo:ADMM} (ADMM) and \Cref{algo:ALADIN} (ALADIN) based on the findings in~\cite{houska2016augmented,engelmann2018toward,jiang2020distributed}.
    Here, CE = central entity.
    }
    \begin{tabular}{llll}
         & ADMM & ALADIN & ALADIN + BFGS  \\ \hline
        Objective function & convex & non-convex & non-convex \\
        Constraints & convex & non-convex & non-convex \\
        Convergence guarantee & no & yes &  yes \\
        Convergence rate & (sub-) linear & quadratic & superlinear \\
        Information local~$\to$~CE & $\sum_{i \in \cS} n_i $ & $\sum_{i \in \cS} \tfrac{n_i(2n_i +3)}{2}$ & $\sum_{i \in \cS} \tfrac{n_i(n_i +4)}{2}$ \\
        Information CE~$\to$~local & $\sum_{i \in \cS} n_i $ & $\sum_{i \in \cS} 2n_i$ & $\sum_{i \in \cS} 2 n_i$ \\
        Applied to energy systems  & \cite{erseghe2014distributed,braun2018hierarchical,mbuwir2020distributed} & \cite{engelmann2018toward,jiang2020distributed} & \cite{engelmann2018toward,zhai2022distributed} \\ \hline
    \end{tabular}

    \label{Tab:CompareADMMandALADIN}
\end{table}
The communication effort is derived in the context of optimal power flow, where~$\cS$ is the set of all sub-networks, cf.~\Cref{Sec:OPF} and~\cite[Tab.~II]{engelmann2018toward}.
The difference in communication effort is mainly due to the second order information used in ALADIN.
While ADMM only communicates local minimizers in both directions (local~$\to$~CE and CE~$\to$~local), the ALADIN algorithm also has to share information of second order derivatives (local~$\to$~CE), and minimizers of~\eqref{eq:dQP} (CE~$\to$~local).

\ \\
Due to its flexibility and convergence properties, \Cref{algo:ALADIN} is widely considered in the literature, see for instance the recent works~\cite{muhlpfordt2021distributed,chanfreut2023aladin,dai2023hypergraph} and~\cite{Jiang22} for applications in the energy management context, where also a comparison between ADMM and ALADIN is presented.
Moreover, in~\cite{jiang2017distributed,shi2018distributed} ALADIN is utilized for vehicle coordination at traffic
intersections; to name but two exemplary instances.

 A potential bottleneck of ALADIN (\Cref{algo:ALADIN}) is its  coupling QP. It involves all decision variables, an active set detection, and a central coordinator. Hence  the original ALADIN variant presented in \Cref{algo:ALADIN} has seen a number of modifications. 
  A trick to overcome the issue of active set detection is to formulate the coupling QP with all inequality constraints, see \cite{botkin2021distributed} for a case study on EV charging. This latter trick, however, does not resolve the scaling issue. 

 To avoid the need for a central coordinator, one can exploit the observation that the coupling QP has the same partially separable structure as the original problem. In \cite{engelmann2020decomposition} a second inner layer of problem decomposition is applied to the convex coupling QP. On this inner level, one may apply ADMM or distributed variants of the conjugate gradient method \cite{engelmann2021essentially}. Moreover, one can project the coupling QP onto the consensus constraints, which ensures that the inner problem scales in the number of coupling variables. This later aspect is of interest for energy grids, which are typically coupled only in a few variables describing line flows or virtual nodes.  We refer to \cite{engelmann2022distributed} for details.  

 Similar to ADMM, by now several algorithmic variants have been proposed. Indeed there is a trend of ALADIN becoming a family of algorithms, cf., e.g.,
\cite{houska2016augmented,engelmann2022aladin,burk2020distributed,botkin2021distributed,houska2021distributed}.

\section{Applications} \label{Sec:Applications}
Next, we present and discuss two prominent applications of the previously introduced optimization methods, 
namely Nonlinear Model Predictive Control of energy distribution in \Cref{Sec:MPC}, and
Optimal Power Flow in \Cref{Sec:OPF}.
While the latter can be considered as steady-state problem, the first application explicitly considers dynamical components.

\subsection{Nonlinear Model Predictive Energy Management} \label{Sec:MPC}
Although invented decades ago, Model Predictive Control (MPC) nowadays is still a topic of ongoing research.
The underlying idea in MPC is to iteratively solve optimal control problems to determine the input. 
To this end, a model to predict the near future behavior of the dynamical system in consideration is required to find controls with minimal costs with respect to a given cost function on the prediction horizon, see for example~\cite{grune2017nonlinear} for an introduction.
Due to its capability to take state and input constraints into account, MPC has been successfully applied in many applications, cf.~\cite{allgower2004nonlinear,vazquez2014model,sultana2017review,hans2018hierarchical,schwenzer2021review}.

In the context of power grids, a typical objective is to avoid rapid changes in the power demand, i.e., the grid operator aims at a flat aggregated power demand profile.
However, since more and more renewable energy sources become part of the overall power network, the overall grid cannot be modeled accurately, and hence, the idea of MPC is not directly applicable.
In particular, due to many small-scale power generation devices, e.g., private PV, the generation and demand profile strongly depends on the local conditions, e.g., sun and wind. Here, recent progress in time-series prediction may be used to generate day-ahead predictions, see, e.g., \cite{sharadga2020time,richter2021day,viehweg2023parameterizing}.
The idea is to consider small-scale semi-autonomous subsystems of the overall power net, and to control the generation-demand balance based on the local sources and conditions, see for example \cite{worthmann2015distributed,braun2016distributed,sauerteig2020towards,braun2018hierarchical,
jiang2020distributed,braun2023model} and references therein, respectively. 

\ \\
In virtue of~\cite{jiang2020distributed}, we consider the following situation as an example.
A set of~$N \in \N$ residential energy systems (private or small-scale commercial) are coupled via a grid operator.
The grid operator (the central entity~CE) has to balance demand and surplus of energy in the smart grid.
A particular component in the smart grid are energy storage devices, e.g., batteries.
In~\cite{jiang2020distributed} the following charging and discharging dynamics for a single energy storage device~$i \in \{1,\ldots,N\}$ are stated
\begin{subequations} \label{eq:MPCBatteryDynamics}
    \begin{align}
        x_i(t+1) &= \alpha_i x(t) + \tau (\beta_i u_i^+(t) + u_i^-(t)), \label{eq:MPCBatteryState} \\
        z(t) &= w_i(t) + u_i^+(t) -\gamma_i u_i^-(t), \label{eq:MPCBatteryDemand}
    \end{align}
\end{subequations}
where~$t \in \N$ is the current time instant, and~$\alpha,\beta,\gamma \in (0,1]$ model efficiency of self-discharge and conversion. 
In~\eqref{eq:MPCBatteryState} the state~$x_i(t) \in \R$ is the state of charge of the~$i^{\rm th}$ battery at time~$t \in \N$, $z_i(t) \in \R$ is the power demand in~\eqref{eq:MPCBatteryDemand}, $w_i(t) \in \R$ is the net energy consumption disregarding the battery, and~$\tau>0$ is the sampling time.
The state of charge~$x_i$ is influenced by the charging~$u_i^+$ and discharging~$u_i^-$ rate, which is considered to be independently controllable, respectively. Hence, we may control the battery by
\begin{equation*}
    u_i(t) = \begin{pmatrix}
        u_i^+(t) \\u_i^-(t)
    \end{pmatrix} \in \R^2.
\end{equation*}
The states~$x_i$, and the charging/discharging rate are subject to constraints.
For battery capacity~$C_i \ge 0$, and bounds on the charging~$\bar u_i \ge 0$ and discharging rate~$\underline u_i \le 0$
typical constraints read, for $i=1,\ldots,N$,
\begin{equation}  \label{eq:MPC_constraints}
    \begin{aligned}
                0 & \le&& x_i(t) &&\le  C_i, \\
                0 &\le &&u_i^+(t) &&\le  \bar u_i, \\
    \underline u_i & \le & &u_i^-(t) &&\le  0, \\
    0 & \le & \frac{u_i^-(t)}{\underline u_i} &+ \frac{u_i^+(t)}{\bar u_i} &&\le  1.
    \end{aligned}
\end{equation}
Invoking linearity in~\eqref{eq:MPCBatteryDynamics}, the state~$x_i$ at time instant~$t+k$ reads
\begin{equation*}
    x_i(t+k) = \alpha_i^{t} x_i(k) + \tau \sum_{j=0}^{t-1} \alpha_i^{t-i-j} \begin{bmatrix} \beta_i & 1 \end{bmatrix} u_i(j).
\end{equation*}
Hence, on the prediction horizon~$0<T \in \N$, for appropriate expressions $D_i \in \R^{8T \times 2 T}$ and $d_i \in \R^{8T}$, the constraints on both, the state and the control~\eqref{eq:MPC_constraints}, can be compactly written as
\begin{equation*}
    D_i u_i \le d_i, \quad i \in \{1,\ldots,N\},
\end{equation*}
which in particular means that on the prediction horizon~$T$ we have polyhedral constraints on the states and controls.
The overall net consumption is given by
\begin{equation*}
    W(t) = \frac{1}{T} \sum_{j = t-T+1}^t \sum_{i=1}^N w_i(j)
\end{equation*}
for~$t \in \{ k,\ldots k+T-1 \}$ with current time step~$k \ge T-1$.
Setting~$ \bar z(t) = \sum_{i=1}^N z_i(t))$, in virtue of MPC, for~$\sigma_0 > 0$ we may formulate the global consensus objective function~$k_0: \R^T \to \Rp$
\begin{subequations} \label{eq:MPC_costs}
\begin{equation} \label{eq:MPC_CE_costs}
    k_0(\bar z) := \frac{\sigma_0}{T N^2} \sum_{j=k}^{k+T-1} (\bar z(j) - W(j))^2 = \frac{\sigma_0}{T N^2} \| \bar z  - W \|^2,
\end{equation}
which describes the deviation between the aggregated power demand and the overall net consumption on the time horizon~$T$.
In the following, we use the predicted demand $z_i := (z_i(t),\ldots,z_i(t+T-1))^\top \in \R^{T}$ 
and the expected net consumption $w_i := (w_i(t),\ldots,w_i(t+T-1))^\top \in \R^{T}$ on the prediction horizon~$T$.

To address local balance between consumption~$w_i$ and demand~$z_i$, and to account for charging/discharging costs, we introduce the local objective functions~$k_i : \R^T \times \R^{2 T} \to \Rp$, $i \in \{1,\ldots,N\}$ as 
\begin{equation} \label{eq:MPC_local_costs}
    k_i(u_i) := \frac{\sigma_i}{2} \big( \| z_i - w_i\|^2 + \| u_i\|^2 \big) = \frac{\sigma_i}{2} \| u_i \|^2_{Q_i}, 
\end{equation}
\end{subequations}
where we used~\eqref{eq:MPCBatteryDemand} to rewrite~$z_i - w_i$. We set 
    $Q_i := \sigma_i  I_{2 T} + \sigma_i I_{T} \otimes \left[\begin{smallmatrix} 1 & \gamma_i \\ \gamma_i & \gamma_i^2 \end{smallmatrix}\right]$.
The factor~$\sigma_i$ penalizes the use of the~$i^{\rm th}$ battery.
With the costs introduced in~\eqref{eq:MPC_costs} and setting $\bar w := \sum_{i=1}^N w_i$ we formulate the optimal control problem
\begin{equation} \label{eq:MPC_OCP}
    \begin{aligned}
         \MIN{\bar z, u} & k_0 (\bar z) + \sum_{i=1}^N k_i(u_i)  \\
          \text{subject to}  &
          \\  \bar z &= \bar w + \sum_{i=1}^N A_i u_i,  & & | \, \lambda_i \\
        d_i & \ge  D_i u_i, && i \in \{1,\ldots,N \} .
    \end{aligned}
\end{equation}
To efficiently solve the optimal control problem~\eqref{eq:MPC_OCP}, in~\cite{jiang2020distributed} ALADIN is utilized.
In each of the~$N$ residential energy systems the following is solved locally and in parallel
\begin{subequations}
\begin{equation} 
    \MIN{v_i} k_i(v_i) - \langle A_i v_i, \lambda \rangle + \frac{1}{2} \| v_i - u_i\|^2_{Q_i} \ \ \text{ subject to } d_i \ge \ D_i v_i \,| \, \kappa_i,
\end{equation}
according to \Cref{algo:ALADIN} Step~1,
where we use the current primal~$u_i$ and dual~$\lambda$ iterate.
First-order optimality $\nabla k_i(v_i) + \langle D_i, \kappa_i \rangle - \langle A_i, \lambda \rangle + Q_i(v_i - u_i) = 0$ allows evaluation of the modified gradient as $g_i = \langle A_i, \lambda \rangle - Q_i(v_i - u_i) $.
The local solutions~$v_i$ are communicated to the CE, where the equality constrained QP
\begin{equation} \label{eq:MPC_ALADIN_QO}
    \begin{aligned}
        \MIN{\bar z,u,s} & k_0(\bar z) + \sum_{i=1}^N \frac{1}{2} \| u_i-v_i\|^2_Q + \langle g_i, u_i \rangle + \frac{\mu_i}{2} \|s_i\|^2 \\
                  \text{ subject to }  & &
          \\  \bar z &= \bar w + \sum_{i=1}^N A_i u_i   ~\quad  | \, \lambda^+ \\
         s_i &= D_i^{\rm act} (u_i - v_i)  \quad | \, \kappa_i^{\rm QP}~ i \in \{1,\ldots,N \} 
    \end{aligned}
\end{equation}
\end{subequations}
is solved.
Here, $D_i^{\rm act} \in \R^{8 T \times 2 T}$ is the Jacobian of the active constraints at local solutions, meaning $D_i^{\rm act} v_i = d_i$ for all $i = 1,\ldots , N$.

For the above problem it is proven in~\cite[Prop.~1]{jiang2020distributed} that the optimization problem~\eqref{eq:MPC_OCP} has a unique solution, and further that for suitable updates of the parameters~$\mu_i$ in~\eqref{eq:MPC_ALADIN_QO} the solution locally converges to the optimum with quadratic rate, cf.~\cite[Prop.~2]{jiang2020distributed}.
Moreover, under regularity assumptions on the solution and a certain Armijo step-size condition \cite[Thm.~1]{jiang2020distributed} ensures global convergence of~\eqref{eq:MPC_OCP} to the optimal solution~$u^\star$.

\subsection{Optimal Power Flow} \label{Sec:OPF}
The notion of Optimal Power Flow (OPF) refers to one of the most important optimization problems in power systems,
where the power grid is described via stationary power flow equations. These equations model routing and distribution as well as loss of power within electricity networks. OPF computations and other optimization problems involving the power flow equations as constraints are applied for a number of different problems in power systems, ranging from the computation of generator setpoints via grid planning to state estimation and others. We refer to \cite{frank2016introduction,faulwasser2018optimal} for tutorial introductions.

One typical variant of OPF considers the computation of setpoints for grid supporting generators, e.g., for so-called redispatch problems when the outcome of  electricity markets is incompatible with physical requirements and operational constraints of the grid, cf.~\cite{faulwasser2018optimal}.
Since stability of the power grid is essential for availability of electrical power, OPF has been and is still a topic of research, see, e.g., \cite{muhlpfordt2017solving,muhlpfordt2019chance} and references therein, respectively.
Introductory overviews and reports on recent developments can be found in the recent work~\cite{bienstock2022mathematical}, and \cite{frank2012optimal,frank2012optimalII,capitanescu2016critical}.
In this section, we consider the rather simple formulation OPF under consideration in~\cite{faulwasser2018optimal}.

\ \\
The stationary behavior of an AC electrical network can be modeled using the triple
$(\cN,\cG,Y)$, where $\cN = \{1,\ldots,N\}$ is the set of nodes (also called buses), the (non-empty) set~$\cG$ describes the
available generators, and $Y = G + j B \in \C^{N \times N}$ is the bus admittance matrix, where $j^2 = -1$.
We introduce OPF for AC systems and then briefly show DC as a special case.

\paragraph{AC power flow}
In a balanced symmetric three-phase AC system every bus $l \in \cN$ can be described by its voltage magnitude~$v_l$, the voltage phase~$\theta_l$, and the net active~$p_l$ and reactive power~$q_l$.
With these quantities we may recall from~\cite{faulwasser2018optimal} the steady-state AC power flow equations for all~$l \in \cN$
\begin{subequations} \label{eq:OPF:PowerFlow}
    \begin{align}
        p_l &= v_l \sum_{k \in \cN} v_k (G_{lk} \cos(\theta_{lk}) + B_{lk} \sin(\theta_{lk}) ), \\
        q_l &= v_l \sum_{k \in \cN} v_k (G_{lk} \sin(\theta_{lk}) - B_{lk} \cos(\theta_{lk}) ), 
    \end{align}
\end{subequations}
where~$\theta_{lk} : = \theta_l - \theta_k$ is the phase angle difference, and $G_{lk},B_{lk}$ denote the admittance of the bus connecting~$l \in \cN$ and $l \neq k \in \cN$, respectively.
Since~$\theta_{lk}$ is defined as differences, we set one bus as reference~$\theta_{l_0}$, w.l.o.g~$1=l_0 \in \cN$.
To simplify notation and exposition we assume $\cG \subseteq \cN$.
For each bus~$l \in \N$ the net apparent power is given by
\begin{equation*}
    s_l := p_l + j q_l = \begin{cases}
        (p_l^g - p_l^d) + j(q_l^g - q_l^d), & l \in \cG, \\
        -(p_l^d + j q_l^d), & \text{else},
    \end{cases}
\end{equation*}
where the superscript~``d'' indicates demand, and~``g'' refers to generation.
For generator nodes~$l \in \cG$ the generated power injections~$p_l^g, q_l^g$ are the control variables,
and for $m \in \cN$ the power sinks or sources~$p_m^d, q_m^d$ are uncontrollable.

\paragraph{DC power flow}
Assuming a lossless line $r_{lk} \approx 0$ for the Ohmic resistance, a constant voltage magnitude~$v_l \approx 1$, 
and small phase differences~$\theta_{lk} \approx 0$ for all~$l,k \in \cN$, we obtain the so-called DC power flow equations as a linear approximation of~\eqref{eq:OPF:PowerFlow} 
\begin{equation}  \label{eq:OPF:PowerFlowLinear}
    p_l = \sum_{k \in \cN \setminus \{l\}} B_{lk}(\theta_l - \theta_k) \ \iff \ p = -B\theta,
\end{equation}
where $B = Im(Y)$, and no conductance is present.

\paragraph{Optimal power flow}
To formulate an optimization problem in a usual way,
we introduce the state $x \in \R^{n_x}$, the disturbance $d \in \R^{n_d}$, and the control input~$u \in \R^{n_u}$ for AC and DC configuration, respectively, by
\begin{subequations} \label{eq:OPF:StatesDisturbancesControls}
    \begin{align}
        x^{\rm AC} &:= (v_l, \theta_l)_{l \in \cN} \in \R^{n_x^{\rm AC}}, \ 
        d^{\rm AC} := (p_l^d, q_l^d)_{l \in \cN} \in \R^{n_d^{\rm AC}},  \
        u^{\rm AC} := (p_l^g, q_l^g)_{l \in \cG} \in \R^{n_u^{\rm AC}}, \\
        x^{\rm DC} &:= (\theta_l)_{l \in \cN} \in \R^{n_x^{\rm DC}}, \ 
        d^{\rm DC} := (p_l^d)_{l \in \cN} \in \R^{n_d^{\rm DC}},  \
        u^{\rm DC} := (p_l^g)_{l \in \cG} \in \R^{n_u^{\rm DC}}. 
    \end{align}
\end{subequations}
In the following $\sigma \in \{\rm AC, DC\}$ indicates the configuration.
Note that the dimensions differ for AC and DC cases as sketched above are as follows
\begin{equation*} \label{eq:OPF:Dimensions}
\begin{array}{l|l|l}
     \sigma & \rm AC & \rm DC  \\ \hline
     n_x^\sigma & 2 |\cN| & |\cN| \\
     n_d^\sigma & 2 |\cN| &  |\cN| \\
     n_u^\sigma & 2 |\cG| & |\cG|
\end{array}
\end{equation*}
Using the triple $(x,d,u)$ from~\eqref{eq:OPF:StatesDisturbancesControls} we may write the power flow equations~\eqref{eq:OPF:PowerFlow} compactly as
\begin{equation} \label{eq:OPF:AbstractPowerFlow}
        F^\sigma : \R^{n_x^\sigma} \times \R^{n_d^\sigma} \times \R^{n_u^\sigma} \to \R^{N^\sigma}, \quad  F^\sigma(x,u;d) = 0, 
\end{equation}
where for AC we have a nonlinear map, and for DC it simplifies to
\begin{equation*}
    F^{\rm DC}(x^{\rm DC},u^{\rm DC};d^{\rm DC}) = u^{\rm DC} - d^{\rm DC} + Bx^{\rm DC} = 0.
\end{equation*}
Note that in this case the dimension of the controls coincides with the dimension of sinks/sources.
Having the abstract formulation~\eqref{eq:OPF:AbstractPowerFlow} at hand, we introduce the 
power flow manifold 
\begin{equation*}  
    \cF^\sigma(d) := \setdef{(x, u) \in \R^{n_x^\sigma} \times \R^{n_u^\sigma}}{ F^\sigma(x,u;d) = 0},
\end{equation*}
which represents all possible solutions of the power flow equations~\eqref{eq:OPF:PowerFlow}, respectively~\eqref{eq:OPF:PowerFlowLinear} for a given disturbance~$d$, i.e., 
the solutions~$(x,u)$ are parameterized by~$d$.
In the following, we assume that the input~$u$ belongs to the box constraints set

   \begin{equation*} 
        \cU^\sigma := \setdef{(u_l^\sigma)_{l \in \cG} \in \R^{n_u^\sigma}}{ 
        \begin{cases}
            p_l^g \in [\underline{p}_l^g, \bar p_l^g], \ q_l^g \in [\underline{q}_l^g, \bar q_l^g] \ \forall \, l \in \cG, & \sigma = \rm AC \\
            p_l^g \in [\underline{p}_l^g, \bar p_l^g],  \ \forall \, l \in \cG, & \sigma = \rm DC \\
        \end{cases}  
        },
   \end{equation*} 
   and the state~$x$ is restricted to the set
   \begin{equation*} 
         \cX := \setdef{(x_l^\sigma)_{l \in \cN} \in \R^{n_x^\sigma} }{ 
         \begin{cases}
         v_l \in [\underline{v}_l, \bar v_l], \ \theta_{1} = 0\ \ \forall \, l \in \cN , & \sigma = \rm AC  \\
          \theta_{1} = 0 \ \ \forall \, l \in \cN , & \sigma = \rm DC  
         \end{cases}
         }.      
   \end{equation*}
   In real applications, there will be further restrictions on the line flows, which are typically stated as constraints on the power
   \begin{equation*} 
       | s_{lk} | = \sqrt{p_{lk}^2 + q_{lk}^2} \le | \bar s_{lk}| \ \ \forall \, (l,k) \in \cL,
   \end{equation*}
for AC, where the set~$\cL \subseteq \cN \times \cN$ represents all connecting lines, and
the active~$p_{lk}$ and reactive power~$q_{lk}$ across the line from~$l \in \cN$ to $l \neq k \in \cN$, i.e.,$(l,k) \in \cL$, is
    \begin{align*}
        p_{lk} &= \phantom{-} v_l^2 g_{lk} - v_l v_m \big( g_{lk} \cos(\theta_{lk}) + b_{lk} \sin(\theta_{lk}) \big), \\
        q_{lk} &= -v_l^2 b_{lk} - v_l v_m \big( b_{lk} \cos(\theta_{lk}) - g_{lk} \sin(\theta_{lk}) \big).
    \end{align*}
Since the constraints on line power flow  depend on the state only, we formulate these compactly 
\begin{equation*}  
\cC^\sigma := \setdef{ x \in \R^{n_x^\sigma} }{ 
\begin{cases}
\sqrt{p_{lk}^2 + q_{lk}^2} \le | \bar s_{lk}| \ \ \forall \, (l,k) \in \cL, & \sigma = \rm AC, \\
- \diag((b_{lk})_{(l,k) \in \cL}) A x \le  (\bar p_{lk})_{(l,k) \in \cL}^\top , & \sigma = \rm DC
\end{cases}
},
\end{equation*}
where~$A \in \R^{|\cL| \times N}$ is the graph incidence matrix in case of DC.
Note that for DC the constraints~$\cC^{\rm DC}$ are obtained from Kirchhoff's laws,
and we have $\cC^{\rm DC} \subset \R^{n_x^{\rm DC}}$.

\ \\
With the notation above, we arrive at the optimization problem
\begin{equation}  \label{eq:OPF_OPT_abstract}
    \begin{aligned}
        \MIN{(x^\sigma,u^\sigma) \in \R^{n_x^\sigma+n_u^\sigma}} \ &J(u^\sigma) \\
        \text{subject to} \quad  (x^\sigma, u^\sigma) & \in \cF^\sigma(d), \\
        u^\sigma & \in \cU^\sigma, \\
        x^\sigma & \in \cX^\sigma \cap \cC^\sigma.
    \end{aligned}
\end{equation}
\paragraph{AC-OPF}
Common control objectives in the AC setting are to reduce transmission losses,
to avoid large injections of reactive power, and to minimize costs of active power generation.
The latter is under consideration in~\cite{engelmann2018toward}, which explicitly results in the following AC optimal power flow problem (AC-OPF)
\begin{subequations} \label{eq:OPF_OPT}
    \begin{align}
    \MIN{\theta,v,p,q} &\sum_{i \in \cG} \alpha_i p_i^2 + \beta_i p_i + \gamma_i, \label{eq:OPF_costs_objective} \\
    \text{subject to}~\eqref{eq:OPF:PowerFlow} & \  \text{and for all} \ k \in \cN \nonumber \\
     &\left.  \begin{array}{l}
         \underline p_i  \le p_i \le \bar p_i, \ \forall i \in \cG, \\
    \underline q_i  \le q_i \le \bar q_i, \ \forall i \in \cG, \\
    \underline v_i  \le v_i \le \bar v_i, \ \forall i \in \cN, \\
    \end{array} \right\} \label{eq:OPF_costs_constraints} \\
    & \ \ \ v_1  = 1 , \ \theta_1 = 0, \label{eq:OPF_costs_reference}
    \end{align}
\end{subequations}
where~$\alpha_i, \beta_i, \gamma_i > 0$.
The terms in~\eqref{eq:OPF_costs_objective} penalize generation of power, the constraints~\eqref{eq:OPF_costs_constraints} are explicit versions of~$u \in \cU^{\rm AC}$ and~$x \in \cX^{\rm AC} \cap \cC^{\rm AC}$ in~\eqref{eq:OPF_OPT_abstract}, and~\eqref{eq:OPF_costs_reference} defines the reference node.
Note that~\eqref{eq:OPF_OPT} contains all nodes, all lines and all constraints in one optimization problem.

In real applications it is often reasonable to consider geographically separated areas as semi-autonomous grids,
in the sense that they are only connected to other areas (grids) by single lines.
In this case the power flow can be considered to be controlled within each grid separately, respectively; while the connecting lines encode constraints on some buses.
With this idea, we may formulate the optimal power flow problem as a task of distributed optimization.

Aiming at representing the AC-OPF~\eqref{eq:OPF_OPT} in the form of~\eqref{eq:SeparableCoupledProblem}, we partition the set of all buses~$\cN$ into~$S$ subsets with~$\cS = \{1,\ldots,S\}$ such that we have the set $\cN_i = \{ n_i^1,\ldots,n_i^{N_i} \} \subset \cN$, where $\cN = \cup_{i \in \cS} = \cN$ and $\cN_i \cap \cN_j = \emptyset$ for $i \neq j$.
The last condition prohibits an overlap of the sub-networks.
To nevertheless account for the connections between the sub-networks, and thus to obtain an equivalent representation of~\eqref{eq:OPF_OPT}, we introduce the set~$\cA$ of auxiliary nodes in those lines which connect two neighboring sub-networks.
Doing so, we obtain an enlarged set $\cN^A := \cN \cup \cA$, containing~$A \in \N$ extra nodes.
Following~\cite{engelmann2018toward}, recalling~$x = (v,\theta)$ and~$u=(p,q)$, we arrive at the separated, affinely coupled optimization problem
\begin{subequations}  \label{eq:OPF_FormulationForALADIN}
    \begin{align}
                 \MIN{(x,u)} &  \sum_{i \in \cS} k_i(u_i) \label{eq:OPF_OPT_separated} \\
            \text{subject to}  &   \nonumber    \\  
           & F^{\rm AC}_i(x_i,u_i,d_i) = 0 && | \, \kappa_i \quad \forall \, i \in \cS, \label{eq:OPF_FlowEq} \\
           & \underline x_i \le h(x_i,u_i,d_i) \le \bar x_i && | \, \xi_i \quad \forall \, i \in \cS, \label{eq:OPF_constraints} \\
           & \sum_{i \in \cS} A_i (x_i,u_i,d_i)^\top = 0 && |\,\lambda.  \label{eq:OPF_coupling} 
    \end{align}
\end{subequations}
In the above formulation, \eqref{eq:OPF_OPT_separated} represent the costs of local power generation with
$k_i(u_i) = \alpha_i p_i^2 + \beta_i p_i + \gamma_i$ for $i \in \cG_i = \cN_i \cap \cG$ (set of local generators), 
\eqref{eq:OPF_FlowEq} are the power flow equations~\eqref{eq:OPF:PowerFlow},
constraints~\eqref{eq:OPF_constraints} compactly encode~\eqref{eq:OPF_costs_constraints},
and \eqref{eq:OPF_coupling} establish the consensus coupling constraints
\begin{equation*}
    \forall\, (k,l) \in \cA \times \cA, \ k \neq l\, : \ \theta_k = \theta_l, \quad v_k = v_l, \quad p_k = -p_l, \quad q_k = -q_l.
\end{equation*}
As discussed in~\cite{engelmann2018toward}, problem~\eqref{eq:OPF_FormulationForALADIN} is amenable to \Cref{algo:ALADIN}, cf.~\cite[Alg.~1]{engelmann2018toward}.
It has been proven in~\cite[Thm.~1]{engelmann2018toward} that, under some (technical) assumptions, the application of \Cref{algo:ALADIN} to problem~\eqref{eq:OPF_FormulationForALADIN} is successful in the sense that it terminates after finite iterations for a user defined accuracy tolerance.
Moreover, under assumptions on the step-size update, the result~\cite[Thm.~2]{engelmann2018toward} ensures local quadratic convergence of the iterates to the optimal value.

\paragraph{Multi-stage AC-OPF}
The power flow equations~\eqref{eq:OPF:PowerFlow} are formulated at steady-state, at least from the viewpoint of control theory.
Hence,~\eqref{eq:OPF_OPT_abstract},~\eqref{eq:OPF_OPT},~\eqref{eq:OPF_FormulationForALADIN} are static optimization problems.
In this setting, however, it is not taken into account that a large-scale power plant cannot change its setpoint arbitrarily fast.
To account for these limitations, in~\cite[Sec.~3.1]{faulwasser2018optimal} multi-stage and predictive optimal power flow is under consideration.
We denote with $u_l(t) = (p_l^g(t), q_l^g(t))$ the control at time instant~$t \in \N$. 
Then, the dynamical limitations (typically only on the active power generation) can, e.g., be formulated as $p^g_l(t+1) - p^g_l(t) \in [\Delta \underline p_l, \Delta \bar p_l]$ for~$l \in \cG$,
leading to control constraints 
\[
u(t+1) - u(t) = \delta u(t) \in \delta \cU := (\times_{l \in \cG} [\Delta \underline p_l, \Delta \bar p_l]) \times \R^{n_u^{\sigma} - |\cG|} \subset \R^{n_u^\sigma}.
\]
Exemplary, for a time horizon~$0<\hat T \in \N $ and $T = \{ 0,\ldots, \hat T \}$, in virtue of~\eqref{eq:OPF_OPT_abstract} we may then state the multi-state AC-OPF
\begin{equation*}
    \begin{aligned}
        \MIN{(x(\cdot),u(\cdot),\delta u(\cdot)) \in \R^{(n_x^\sigma+2 n_u^\sigma) \hat T}} \ & \sum_{k \in T} J(u(k)) + \| \delta u(k)\|^2 \\
        \text{subject to} \ \forall \, k \in T & \\
        u(k+1) &= u(k) + \delta u(k), \\
        \delta u(k) &\in \delta \cU \\ 
         (x(k), u(k)) & \in (\cX^{\rm AC} \times \cU^{\rm AC}) \cap \cF^{\rm AC}(d(k)), \\
        x(k) & \in \cC^{\rm AC},
    \end{aligned}
\end{equation*}
where the second addend in the objective function penalizes frequency and amount of power generation.
The multi-state AC-OPF already has a flavor of predictive control presented in \Cref{Sec:MPC}, and its form is similar to that of an optimal control problem, cf.~\cite{faulwasser2020toward}.

\begin{figure}
\begin{subfigure}[t]{0.50\linewidth}
\includegraphics[width=\textwidth]{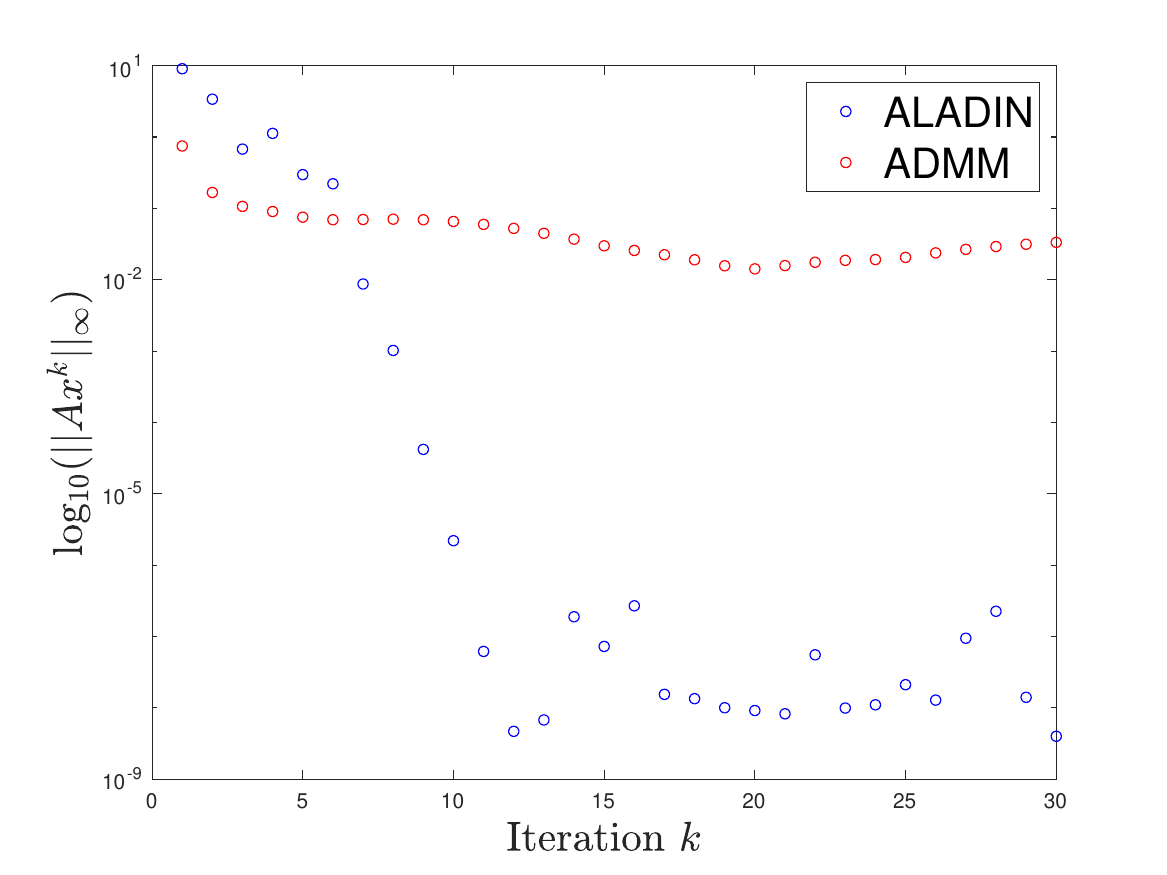}
\subcaption{Constraint violation.}
\label{Fig:CompALADIN_ADMM}
    \end{subfigure}
    \begin{subfigure}[t]{0.33\linewidth}
    \includegraphics[width=\textwidth, angle=90]{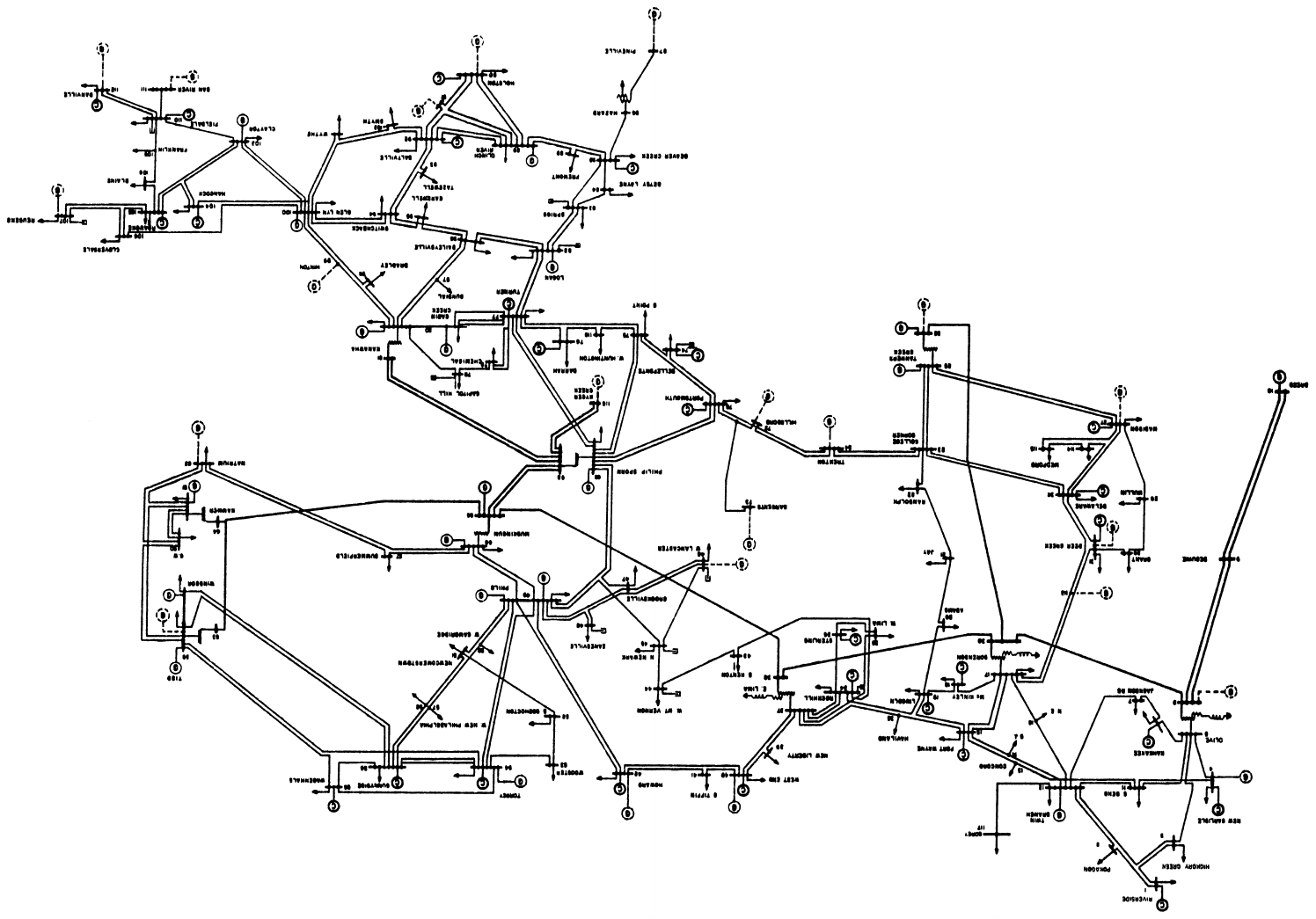}
\subcaption{IEEE 118 bus test system.}
\label{Fig:IEEE_TestSystem}
\end{subfigure}
\caption{Left: logarithmic constraint violation for ADMM and ALADIN. Right: example of the IEEE 118 bus test system.}
\label{Fig:NumExample}
\end{figure}

\paragraph{Numerical example}
To showcase functioning of the proposed methods, we briefly present a numerical example of an AC-OPF problem~\eqref{eq:OPF_FormulationForALADIN}.
We simulate the IEEE 118-bus test case, considered in~\cite{engelmann2018toward}.
\Cref{Fig:IEEE_TestSystem} shows the map of the IEEE 118-bus test system.
We compare application of \Cref{algo:ADMM} (ADMM) with \Cref{algo:ALADIN} (ALADIN).
Exemplary, in \Cref{Fig:CompALADIN_ADMM} violation of the coupling constraints~\eqref{eq:OPF_coupling} is depicted.
It can be seen that \Cref{algo:ALADIN} is superior to \Cref{algo:ADMM} in terms of constraint satisfaction.
The simulation has been performed using the \textsc{Matlab} toolbox \textsc{ALADIN$-\alpha$}~\cite{engelmann2022aladin}.
We refer to~\cite{engelmann2018toward,engelmann2022aladin}, where more detailed and comprehensive numerical case studies are conducted.
Therein, also smaller (smallest 5-bus system) and larger (largest 300-bus system) are considered, and different variants of ALADIN are compared with respect to number of iterations against relative error, communication effort, and convergence rate.
Notice that the \textsc{ALADIN$-\alpha$} toolbox~\cite{engelmann2022aladin} with many instructive examples for different distributed optimization problems is available open-source~\cite{ALADIN-toolbox}.

\section{Outlook and open problems} \label{Sec:Outlook}
ADMM and ALADIN are only two distributed algorithms that can be applied to power and energy systems. Indeed, more recent advances have considered distributing interior point methods and sequential quadratic programming, see \cite{engelmann2021essentially} and \cite{stomberg2022decentralized} for details and for numerical results for AC-OPF problems (also in comparison to ADMM and ALADIN).
A very much open problem in distributed optimization is the formal analysis of the interplay between the actually considered decomposition and the realized convergence speed of the algorithms. 
In general, ALADIN (\Cref{algo:ALADIN}) exhibits convergence guarantees for non-convex optimization problems and the guarantees often rely on a globalization step which as such is difficult to distribute between agents.

Some open-source codes and tools for distributed optimization are available. In case of ALADIN~\cite{engelmann2022aladin} provides an open-source \textsc{MATLAB} implementation of different variants tailored for algorithmic prototyping, while \cite{muhlpfordt2021distributed} provides code tailored to power system applications and TSO-DSO coordination. 
Similarly, ADMM codes have been published, cf. \cite{engelmann2022aladin,stomberg2023cooperative} and references therein. However, in comparison to the quite large number of powerful open-source codes for solving NLPs and optimal control problems, there is still demand for further tool development.

For both---numerical optimization algorithms and for modeling and control of distributed energy systems---data-driven and learning approaches are of increasing interest~\cite{bertozzi2023application}.
With respect to optimization, see, e.g., \cite{baumann2019surrogate}, where recurrent neural networks are employed to reduce the number of iterations and \cite{zeng2022reinforcement} for learning of parameters for ADMM. 
for recent case studies on power systems. Within this framework, recent extensions to stochastic systems are key to properly take uncertainties into account~\cite{faulwasser2023behavioral}.

From an applications perspective, scalable algorithms for distributed stochastic OPF, distributed approaches to reactive power dispatch (which involves discrete decision variables~\cite{murray2018hierarchical}), and the solution of problems involving coupled power systems and gas networks are of interest. Moreover, efficient numerical implementations which can be applied on small-scale computational units (a.k.a. embedded systems) in distribution grids and multi-energy problems~\cite{sass2020model} are of growing relevance. Another topic of interest is to go from separate open-source grid models~\cite{meinecke2020simbench} and open-source datasets of load and generation time series such as \cite{ratnam2017residential,spalthoff2019simbench} to  coupled grid models and time-series data, which  appear to be not yet available. Likewise, open-source datasets and benchmarks for multi-energy systems such as~\cite{sass2020model} appear to be still rare.

\paragraph{Acknowledgments}
L.~Lanza gratefully acknowledges  support by the Carl Zeiss Foundation (VerneDCt -- Project-ID 2011640173).

\end{document}